\documentclass[a4paper,twoside,11pt]{article}
\usepackage{bbm}
\usepackage{amsfonts}
\usepackage{amssymb}\usepackage{amsmath}
 \numberwithin{equation}{section}
\begin{document}

\title{On the Cauchy problem of a periodic 2-component $\mu$-Hunter-Saxton system}

\author{Jingjing Liu\footnote{e-mail:
jingjing830306@163.com }\\
Department of Mathematics,
Sun Yat-sen University,\\
510275 Guangzhou, China
\bigskip\\
Zhaoyang Yin\footnote{e-mail: mcsyzy@mail.sysu.edu.cn}
\\ Department of Mathematics, Sun Yat-sen University,\\ 510275 Guangzhou, China}
\date{}
\maketitle

\begin{abstract}
In this paper, we study the Cauchy problem of a periodic 2-component
$\mu$-Hunter-Saxton system. We first establish the local
well-posedness for the periodic 2-component $\mu$-Hunter-Saxton
system by Kato's semigroup theory. Then, we derive the precise
blow-up scenario for strong solutions to the system. Moreover, we
present some blow-up results for strong solutions to the system.
Finally, we give a
global existence result to the system.\\

\end{abstract}

\section{Introduction}
\par
Recently, a new 2-component system  was introduced by Zuo in
\cite{d-z} as follows:
\begin{equation}
\left\{\begin{array}{ll}
\mu(u)_{t}-u_{txx}=2\mu(u)u_{x}-2u_{x}u_{xx}-uu_{xxx}+\rho\rho_{x}&-\gamma_{1}u_{xxx},\\
&t > 0,\,x\in \mathbb{R},\\
\rho_{t}=(\rho u)_x+2\gamma_{2}\rho_{x}, &t > 0,\,x\in \mathbb{R},\\
u(0,x) = u_{0}(x),& x\in \mathbb{R}, \\
\rho(0,x) = \rho_{0}(x),&x\in \mathbb{R},\\
u(t,x+1)=u(t,x), & t \geq 0, x\in \mathbb{R},\\
\rho(t,x+1)=\rho(t,x), & t \geq 0, x\in \mathbb{R},\\ \end{array}\right. \\
\end{equation}
where $\mu(u)=\int_{\mathbb{S}}udx$ with
$\mathbb{S}=\mathbb{R}/\mathbb{Z}$ and $\gamma_{i}\in \mathbb{R},$
$i=1,2.$ By integrating both sides of the first equation in the
system (1.1) over the circle $\mathbb{S}=\mathbb{R}/\mathbb{Z}$ and
using the periodicity of $u$, one obtain
$$\mu(u_{t})=\mu(u)_{t}=0.$$
This yields the following periodic 2-component $\mu$-Hunter-Saxton
system:
\begin{equation}
\left\{\begin{array}{ll}
-u_{txx}=2\mu(u)u_{x}-2u_{x}u_{xx}-uu_{xxx}+\rho\rho_{x}&-\gamma_{1}u_{xxx},\\
&t > 0,\,x\in \mathbb{R},\\
\rho_{t}=(\rho u)_x+2\gamma_{2}\rho_{x}, &t > 0,\,x\in \mathbb{R},\\
u(0,x) = u_{0}(x),& x\in \mathbb{R}, \\
\rho(0,x) = \rho_{0}(x),&x\in \mathbb{R},\\
u(t,x+1)=u(t,x), & t \geq 0, x\in \mathbb{R},\\
\rho(t,x+1)=\rho(t,x), & t \geq 0, x\in \mathbb{R},\\ \end{array}\right. \\
\end{equation}
with $\gamma_{i}\in \mathbb{R},$ $i=1,2$. This system is a
2-component generalization of the generalized Hunter-Saxton equation
obtained in \cite{k-l-m}. The author \cite{d-z} shows that this
system is both a bihamiltonian Euler equation and a bivariational
equation.

Obviously, (1.1) is equivalent to (1.2) under the condition
$\mu(u_{t})=\mu(u)_{t}=0.$ In this paper, we will study the system
(1.2) under the assumption $\mu(u_{t})=\mu(u)_{t}=0$.

For $\rho\equiv 0$ and $\gamma=0,$ and replacing $t$ by $-t,$ the
system (1.2) reduces to the generalized Hunter-Saxton equation
(named $\mu$-Hunter-Saxton equation or $\mu$-Camassa-Holm equation)
as follows:
\begin{equation}
-u_{txx}=-2\mu(u)u_{x}+2u_{x}u_{xx}+uu_{xxx},
\end{equation}
which is obtained and studied in \cite{k-l-m}. Moreover, the
periodic $\mu$-Hunter-Saxton equation and the periodic
$\mu$-Degasperis-Procesi equation have also been studied in
\cite{fu, l-m-t} recently. It is worthy to note that the
$\mu$-Hunter-Saxton equation has a very closed relation with the
periodic Hunter-Saxton and Camassa-Holm equations.  For $\mu(u)=0,$
the equation (1.3) reduces to the Hunter-Saxton equation \cite{J-R}
\begin{equation}
u_{txx}+2u_{x}u_{xx}+uu_{xxx}=0,
\end{equation}
modeling the propagation of weakly nonlinear orientation waves in a
massive nematic liquid crystal. Here, $u(t,x)$ describes the
director field of a nematic liquid crystal, $x$ is the space
variable in a reference frame moving with the linearized wave
velocity, $t$ is a slow time variable. The orientation of the
molecules is described by the field of unit vectors $(\cos u(t,x),
\sin u(t,x))$ \cite{Y}. The single-component model also arises in a
different physical context as the high-frequency limit \cite{HHD,
J-Y} of the Camassa-Holm equation, which is a model for shallow
water waves \cite{R-D, J-H} and a re-expression of the geodesic flow
on the diffeomorphism group of the circle \cite{C-K} with a
bi-Hamiltonian structure \cite{F-F} which is completely integrable
\cite{C-M}. The Hunter-Saxton equation also has a bi-Hamiltonian
structure \cite{J-H, P-P} and is completely integrable \cite{B1,
J-Y}. The initial value problem for the Hunter-Saxton equation (1.4)
on the line (nonperiodic case) and on the unit circle
$\mathbb{S}=\mathbb{R}/\mathbb{Z}$ were studied by Hunter and Saxton
in \cite{J-R} using the method of characteristics and by Yin in
\cite{Y} using Kato semigroup method, respectively.

For $\rho\not\equiv0$ , $\gamma_{i}=0, i=1,2$ £¬$\mu(u)=0$ and
replacing $t$ by $-t,$ peakon solutions of the Cauchy problem of the
system (1.2) have been analysed in \cite{C-I}. Moreover, the Cauchy
problem of 2-component periodic Hunter-Saxton system has been
discussed in \cite{l-y, MW}. However, the Cauchy problem of the
system (1.2) has not been studied yet. The aim of this paper is to
establish the local well-posedness for the system (1.2), to derive
the precise blow-up scenario, to prove that the system (1.2) has
global strong solutions and also finite time blow-up solutions.

The paper is organized as follows. In Section 2, we establish the
local well-posedness of the initial value problem associated with
the system (1.2). In Section 3, we derive the precise blow-up
scenario. In Section 4, we present two explosion criteria of strong
solutions to the system (1.2) with general initial data. In Section
5, we give a new global existence result of strong solutions to the
system (1.2).

 \textbf{Notation}  Given a Banach space $Z$, we denote its norm by
 $\|\cdot\|_{Z}$. Since all space of functions are over
 $\mathbb{S}=\mathbb{R}/\mathbb{Z}$, for simplicity, we drop $\mathbb{S}$ in our notations
 if there is no ambiguity. We let $[A,B]$ denote
 the commutator of linear operator $A$ and $B$. For convenience, we
 let $(\cdot|\cdot)_{s\times r}$ and $(\cdot|\cdot)_{s}$ denote the
 inner products of $H^{s}\times H^{r}$, $s,r\in \mathbb{R}_{+}$ and
 $H^{s}$, $s\in \mathbb{R}_{+}$, respectively.

\section{Local well-posedness}
\newtheorem {remark2}{Remark}[section]
\newtheorem{theorem2}{Theorem}[section]
\newtheorem{lemma2}{Lemma}[section]
In this section, we will establish the local well-posedness for the
Cauchy problem of the system (1.2) in $H^{s}\times H^{s-1}$, $s\geq
2$, by applying Kato's theory \cite{Kato 1}.

The condition $\mu(u_{t})=0$ ensures that the first equation in
(1.2) can be recast in the form
$$u_{t}-(u+\gamma_{1})u_{x}=\partial_{x}(\mu-\partial_{x}^{2})^{-1}(2\mu
u+\frac{1}{2}u_{x}^{2}+\frac{1}{2}\rho^{2}),$$ where
$A=\mu-\partial_{x}^{2}$ is an isomorphism between $H^{s}$ and
$H^{s-2}.$ Using this identity, the system (1.2) takes the form of a
quasi-linear evolution equation of hyperbolic type:
\begin{equation}
 \ \ \ \ \ \ \ \  \ \ \ \ \left\{\begin{array}{ll}
u_{t}-(u+\gamma_{1})u_{x}=\partial_{x}(\mu-\partial_{x}^{2})^{-1}&(2\mu
u+\frac{1}{2}u_{x}^{2}+\frac{1}{2}\rho^{2}),\\
&t > 0,\,x\in \mathbb{R},\\
\rho_{t}-(u+2\gamma_{2})\rho_{x}=u_{x}\rho, &t > 0,\,x\in \mathbb{R},\\
u(0,x) = u_{0}(x),&x\in \mathbb{R}, \\
\rho(0,x) = \rho_{0}(x),&x\in \mathbb{R},\\
\rho(t,x+1)=\rho(t,x), & t \geq 0, x\in \mathbb{R},\\
u(t,x+1)=u(t,x), & t \geq 0, x\in \mathbb{R}.\\
\end{array}\right. \\
\end{equation}

Let $z:=\left(
          \begin{array}{c}
            u \\
            \rho \\
          \end{array}
        \right),$ $A(z)=\left(
         \begin{array}{cc}
           -(u+\gamma_{1})\partial_{x} & 0 \\
           0 & -(u+2\gamma_{2})\partial_{x} \\
         \end{array}
       \right)$ and
$$f(z)=\left(
              \begin{array}{c}
                \partial_{x}(\mu-\partial_{x}^{2})^{-1}(2\mu
u+\frac{1}{2}u_{x}^{2}+\frac{1}{2}\rho^{2})\\
                 u_{x}\rho \\
              \end{array}
            \right).$$
Set $Y=H^{s}\times H^{s-1},$ $X=H^{s-1}\times H^{s-2},$
$\Lambda=(\mu-\partial_{x}^{2})^{\frac{1}{2}}$ and $Q=\left(
                                                     \begin{array}{cc}
                                                       \Lambda & 0 \\
                                                       0 & \Lambda \\
                                                     \end{array}
                                                   \right).$
Obviously, $Q$ is an isomorphism of $H^{s}\times H^{s-1}$ onto
$H^{s-1}\times H^{s-2}.$\\

Similar to the proof of Theorem 2.2 in \cite{E-L-Y}, we get the
following conclusion.

\begin{theorem2}
Given $z_{0}=(u_{0},\rho_{0})\in H^{s}\times H^{s-1}$, $s\geq 2,$
then there exists a maximal $T=T(\parallel
z_{0}\parallel_{H^{s}\times H^{s-1}})>0$, and a unique solution
$z=(u, \rho)$ to (2.1) such that
$$ z=z(\cdot,z_{0})\in C([0,T); H^{s}\times H^{s-1})\cap
C^{1}([0,T);H^{s-1}\times H^{s-2}).
$$ Moreover, the solution depends continuously on the initial data, i.e., the
mapping $$z_{0}\rightarrow z(\cdot,z_{0}): H^{s}\times
H^{s-1}\rightarrow C([0,T); H^{s}\times H^{s-1})\cap
C^{1}([0,T);H^{s-1}\times H^{s-2})
$$
is continuous.
\end{theorem2}

Recall that the periodic 2-component Hunter-Saxton system discussed
in \cite{l-y} only has local existence but not local well-posedness
because of the lack of uniqueness. The ambiguity disappears in the
case of the periodic 2-component $\mu$-Hunter-Saxton system from the
Theorem 2.1. This is a very important difference between the
2-component Hunter-Saxton system and the
2-component $\mu$-Hunter-Saxton system.\\

Consequently, we will give another equivalent form of (1.2).
Integrating both sides of the first equation in (1.2) with respect
to $x$, we obtain
$$u_{tx}=-2\mu(u)u+\frac{1}{2}u_{x}^{2}+uu_{xx}-\frac{1}{2}\rho^{2}+\gamma_{1}u_{xx}+a(t),$$
where
$$a(t)=2\mu(u)^{2}+\frac{1}{2}\int_{\mathbb{S}}(u_{x}^{2}+\rho^{2})dx.$$
Using the system (1.2), we have
\begin{align}
&\frac{1}{2}\frac{d}{dt}\int_{\mathbb{S}}(u_{x}^{2}+\rho^{2})dx\\
\nonumber=&\int_{\mathbb{S}}(u_{x}u_{xt}+\rho\rho_{t})dx\\
\nonumber=&-\int_{\mathbb{S}}uu_{txx}dx+\int_{\mathbb{S}}\rho\rho_{t}dx\\
\nonumber=&\int_{\mathbb{S}}2\mu(u)
uu_{x}dx-2\int_{\mathbb{S}}uu_{x}u_{xx}dx-\int_{\mathbb{S}}u^{2}u_{xxx}dx+\int_{\mathbb{S}}u\rho_{x}\rho dx\\
\nonumber&-\gamma_{1}\int_{\mathbb{S}}uu_{xxx}dx+\int_{\mathbb{S}}\rho(u\rho)_{x}dx+2\gamma_{2}\int_{\mathbb{S}}\rho\rho_{x}dx\\
\nonumber=&\int_{\mathbb{S}}u\rho_{x}\rho
dx+\int_{\mathbb{S}}\rho(u\rho)_{x}dx=0.
\end{align}
By $\mu(u)_{t}=\mu(u_{t})=0,$ we have
$$\frac{d}{dt}a(t)=0.$$
For convenience, we let
$$\mu_{0}:=\mu(u_{0})=\mu(u)=\int_{\mathbb{S}}u(t,x)dx,$$
$$\mu_{1}:=\left(\int_{\mathbb{S}}(u_{x}^{2}+\rho^{2})dx\right)^{\frac{1}{2}}
=\left(\int_{\mathbb{S}}(u_{0,x}^{2}+\rho_{0}^{2})dx\right)^{\frac{1}{2}}$$
and write $a:=a(0)$ henceforth. Thus,
\begin{equation}
u_{tx}=-2\mu_{0}u+\frac{1}{2}u_{x}^{2}+uu_{xx}-\frac{1}{2}\rho^{2}+\gamma_{1}u_{xx}+a
\end{equation}
is a valid reformulation of the first equation in (1.2). Integrating
(2.3) with respect to $x$, we get
$$
u_{t}-(u+\gamma_{1})u_{x}=\partial_{x}^{-1}(-2\mu_{0}u-\frac{1}{2}u_{x}^{2}-\frac{1}{2}\rho^{2}+a)+h(t),
$$
where $\partial_{x}^{-1}g(x)=\int_{0}^{x}g(y)dy$ and $h(t):
[0,\infty)\rightarrow \mathbb{R}$ is a continuous function. For the
2-component Hunter-Saxton system, if we follow the above same
procedure, then the arbitrariness of continuous function $h(t)$ will
lead to the non-uniqueness of solution to (1.2). In this paper, the
condition $\mu(u_{t})=0$ implies that $h(t)$ is unique.
Consequently, the solution to (1.2) will be unique.

Thus we get another equivalent form of (1.2)
\begin{equation}
 \ \ \ \ \ \ \ \  \ \ \ \ \left\{\begin{array}{ll}
u_{t}-(u+\gamma_{1})u_{x}=\partial_{x}^{-1}(-2\mu_{0}u-\frac{1}{2}u_{x}^{2}\\-\frac{1}{2}\rho^{2}+a)+h(t),
&t > 0,\,x\in \mathbb{R},\\
\rho_{t}-(u+2\gamma_{2})\rho_{x}=u_{x}\rho, &t > 0,\,x\in \mathbb{R},\\
u(0,x) = u_{0}(x),&x\in \mathbb{R}, \\
\rho(0,x) = \rho_{0}(x),&x\in \mathbb{R},\\
u(t,x+1)=u(t,x), & t \geq 0, x\in \mathbb{R},\\
\rho(t,x+1)=\rho(t,x), & t \geq 0, x\in \mathbb{R},\\
\end{array}\right. \\
\end{equation}
where $\partial_{x}^{-1}g(x)=\int_{0}^{x}g(y)dy$ and $h(t):
[0,\infty)\rightarrow \mathbb{R}$ is a continuous function.

\section{The precise blow-up scenario}
\newtheorem {remark3}{Remark}[section]
\newtheorem{theorem3}{Theorem}[section]
\newtheorem{lemma3}{Lemma}[section]
\newtheorem{definition3}{Definition}[section]
\newtheorem{claim3}{Claim}[section]

In this section, we present the precise blow-up scenario for strong
solutions to the system (1.2).

We first recall the following lemmas.

\begin{lemma3}
\cite{Kato 3} If $r>0$, then $H^{r}\cap L^{\infty}$ is an algebra.
Moreover
$$\parallel fg\parallel_{H^r}\leq c(\parallel  f\parallel_{L^{\infty}}\parallel g\parallel_{H^r}+\parallel
  f\parallel_{H^r}\parallel g\parallel_{L^{\infty}}),
$$
where c is a constant depending only on r.
\end{lemma3}

\begin{lemma3}
\cite{Kato 3} If $r>0$, then
$$\parallel[\Lambda^{r},f]g\parallel_{L^{2}}\leq c(\parallel \partial_{x}f\parallel_{L^{\infty}}
\parallel\Lambda^{r-1}g\parallel_{L^2}+\parallel
\Lambda^r f\parallel_{L^2}\parallel g\parallel_{L^{\infty}}),
$$
where c is a constant depending only on r.
\end{lemma3}

Next we prove the following useful result on global existence of
solutions to (1.2).
\begin{theorem3} Let $z_{0}=\left(\begin{array}{c}
                                u_{0} \\
                                \rho_{0} \\
                              \end{array}
                            \right)
\in H^{s}\times H^{s-1}$, $s\geq 2$, be given and assume that T is
the maximal existence time of the corresponding solution
                  $z=\left(\begin{array}{c}
                                      u \\
                                      \rho \\
                                    \end{array}
                                  \right)$
to (2.4) with the initial data $z_{0}$. If there exists $M>0$ such
that
$$
\| u_{x}(t,\cdot)\|_{L^{\infty}}
+\|\rho(t,\cdot)\|_{_{L^{\infty}}}+\|\rho_{x}(t,\cdot)\|_{_{L^{\infty}}}\leq
M,\ \ t\in[0,T),
$$
then the $H^s\times H^{s-1}$-norm of $z(t,\cdot)$ does not blow up
on [0,T).
 \end{theorem3}
\textbf{Proof}
 Let $z=\left(
                                    \begin{array}{c}
                                      u \\
                                      \rho \\
                                    \end{array}
                                  \right)$ be the solution to
(2.4) with the initial data $z_{0}\in H^s\times H^{s-1},\ s\geq 2$,
and let T be the maximal existence time of the corresponding
solution $z$, which is guaranteed by Theorem 2.1. Throughout this
proof, $c>0$ stands for a generic constant depending only on $s$.

Applying the operator $\Lambda^{s}$ to the first equation in (2.4),
multiplying by $\Lambda^{s} u$, and integrating over $\mathbb{S}$,
we obtain
\begin{equation}
\frac{d}{dt}\|u\|^{2}_{H^{s}}=2(uu_{x},
u)_{s}+2(u,\partial_{x}^{-1}(-2\mu_{0}
u-\frac{1}{2}u_{x}^{2}-\frac{1}{2}\rho^{2}+a)+h(t))_{s}.
\end{equation}
Let us estimate the first term of the right-hand side of (3.1).
\begin{align}
|(uu_{x},u)_{s}|&= |(\Lambda^s(u\partial_{x}u),\Lambda^s u)_{0}|\\
\nonumber&=|([\Lambda^s,u]\partial_{x}u,\Lambda^s
u)_{0}+(u\Lambda^s\partial_{x}u,\Lambda^s u)_{0}|\\
\nonumber&\leq \|[\Lambda^s,u]\partial_{x}u\|_{L^2}\|\Lambda^s
u\|_{L^2}+\frac{1}{2}|(u_{x}\Lambda^s u,\Lambda^s u)_{0}|\\
\nonumber& \leq (c\| u_{x}\|_{L^{\infty}}+\frac{1}{2}\|
u_{x}\|_{L^{\infty}})\|u\|^2_{H^s}\\
\nonumber&\leq c\|u_{x}\|_{L^{\infty}}\|u\|^2_{H^s},
\end{align}
where we used Lemma 3.2 with $r=s$. Let $f\in H^{s-1}, s\geq 2$. We
have
$$
|\partial_x^{-1} f|=|\int_0^{x}f dx|\leq \int_{\mathbb{S}}|f|dx\leq
\|f\|_{L^{2}}
$$
and
$$
 \|\partial_x^{-1} f\|_{L^2}=\left(\int_0^1(\partial_x^{-1}
 f)^2dx\right)^{1/2}\leq
\left(\int_0^1\|f\|_{L^{2}}^2dx\right)^{1/2}=\|f\|_{L^2}.
$$
Thus
$$ \|\partial_x^{-1}f\|_{H^{s}}\leq \|\partial_x^{-1}
f\|_{L^2}+\|f\|_{H^{s-1}}\leq 2\|f\|_{H^{s-1}}.$$ Then, we estimate
the second term of the right-hand side of (3.1) in the following
way:
\begin{align}
&|(\partial_{x}^{-1}(-2\mu_{0}
u-\frac{1}{2}u_{x}^{2}-\frac{1}{2}\rho^{2}+a)+h(t),u)_{s}|\\
\nonumber\leq \ &\|\partial_{x}^{-1}(-2\mu_{0}
u-\frac{1}{2}u_{x}^{2}-\frac{1}{2}\rho^{2}+a)+h(t)\|_{H^{s}}\|u\|_{H^{s}}\\
\nonumber\leq  \ &(\|\partial_{x}^{-1}(-2\mu_{0}
u-\frac{1}{2}u_{x}^{2}-\frac{1}{2}\rho^{2}+a)\|_{H^{s}}
+\|h(t)\|_{H^{s}})\|u\|_{H^{s}}\\
\nonumber\leq \ &(2\|-2\mu_{0}
u-\frac{1}{2}u_{x}^{2}-\frac{1}{2}\rho^{2}+a\|_{H^{s-1}}
+\|h(t)\|_{H^{s}})\|u\|_{H^{s}}\\
\nonumber\leq  \
&(4|\mu_{0}|\|u\|_{H^{s}}+\|u_{x}^{2}\|_{H^{s-1}}+\|\rho^{2}\|_{H^{s-1}}+2\|a\|_{H^{s-1}}+\|h(t)\|_{H^{s}})\|u\|_{H^{s}}\\
\nonumber\leq  \ &c
(\|u\|_{H^{s}}+\|u_{x}\|_{L^{\infty}}\|u_{x}\|_{H^{s-1}}+\|\rho\|_{L^{\infty}}\|\rho\|_{H^{s-1}}+|a|+\max\limits_{t\in[0,T)}|h(t)|)\|u\|_{H^{s}}\\
\nonumber\leq \  &c
(\|u_{x}\|_{L^{\infty}}+\|\rho\|_{L^{\infty}}+1)(\|u\|_{H^{s}}^{2}+\|\rho\|_{H^{s-1}}^{2}+1),
\end{align}
where we used Lemma 3.1 with $r=s-1$. Combining (3.2) and (3.3) with
(3.1), we get
\begin{equation}
\frac{d}{dt}\|u\|^{2}_{H^{s}}\leq
c(\|\rho\|_{L^{\infty}}+\|u_{x}\|_{L^{\infty}}+1)(\|u\|_{H^{s}}^{2}+\|\rho\|_{H^{s-1}}^{2}+1).
\end{equation}
In order to derive a similar estimate for the second component
$\rho$, we apply the operator $\Lambda^{s-1}$ to the second equation
in (2.4), multiply by $\Lambda^{s-1}\rho$, and integrate over
$\mathbb{S}$, to obtain
\begin{equation}
\frac{d}{dt}\|\rho\|_{H^{s-1}}^{2}=2(u\rho_{x},\rho)_{s-1}+2(u_{x}\rho,\rho)_{s-1}.
\end{equation}
Let us estimate the first term of the right hand side of (3.5)
\begin{align*}
&|(u\rho_{x},\rho)_{s-1}|\\
=  \ &|(\Lambda^{s-1}(u\partial_{x}\rho),\Lambda^{s-1}\rho)_{0}|\\
=  \ &|([\Lambda^{s-1},u]\partial_{x}\rho,
\Lambda^{s-1}\rho)_{0}+(u\Lambda^{s-1}\partial_{x}\rho,
\Lambda^{s-1}\rho)_{0}|\\
\leq \
&\|[\Lambda^{s-1},u]\partial_{x}\rho\|_{L^{2}}\|\Lambda^{s-1}\rho\|_{L^{2}}+\frac{1}{2}|(u_{x}\Lambda^{s-1}\rho,\Lambda^{s-1}\rho)_{0}|\\
\leq \
&c(\|u_{x}\|_{L^{\infty}}\|\rho\|_{H^{s-1}}+\|\rho_{x}\|_{L^{\infty}}\|u\|_{H^{s-1}})\|\rho\|_{H^{s-1}}
+\frac{1}{2}\|u_{x}\|_{L^{\infty}}\|\rho\|_{H^{s-1}}^{2}\\
\leq \
&c(\|u_{x}\|_{L^{\infty}}+\|\rho_{x}\|_{L^{\infty}})(\|\rho\|_{H^{s-1}}^{2}+\|u\|_{H^{s}}^{2}),
\end{align*}
here we applied Lemma 3.2 with $r=s-1$. Then we estimate the second
term of the right hand side of (3.5). Based on Lemma 3.1 with
$r=s-1$, we get
\begin{align*}
|(u_{x}\rho,\rho)_{s-1}|\leq \ &\|u_{x}\rho\|_{H^{s-1}}\|\rho\|_{H^{s-1}}\\
\leq \
&c(\|u_{x}\|_{L^{\infty}}\|\rho\|_{H^{s-1}}+\|\rho\|_{L^{\infty}}\|u_{x}\|_{H^{s-1}})\|\rho\|_{H^{s-1}}\\
\leq \
&c(\|u_{x}\|_{L^{\infty}}+\|\rho_{x}\|_{L^{\infty}})(\|\rho\|_{H^{s-1}}^{2}+\|u\|_{H^{s}}^{2}).
\end{align*}
Combining the above two inequalities with (3.5), we get
\begin{equation}
\frac{d}{dt}\|\rho\|^{2}_{H^{s-1}}\leq
c(\|u_{x}\|_{L^{\infty}}+\|\rho\|_{L^{\infty}}+\|\rho_{x}\|_{L^{\infty}})(\|u\|_{H^{s}}^{2}+\|\rho\|_{H^{s-1}}^{2}+1).
\end{equation}
By (3.4) and (3.6), we have
\begin{align*}
&\frac{d}{dt}(\|u\|_{H^{s}}^{2}+\|\rho\|^{2}_{H^{s-1}}+1)\\
\leq \
&c(\|u_{x}\|_{L^{\infty}}+\|\rho\|_{L^{\infty}}+\|\rho_{x}\|_{L^{\infty}}+1)(\|u\|_{H^{s}}^{2}+\|\rho\|_{H^{s-1}}^{2}+1).
\end{align*}
An application of Gronwall's inequality and the assumption of the
theorem yield
$$(\|u\|_{H^{s}}^{2}+\|\rho\|^{2}_{H^{s-1}}+1)\leq\exp(c(M+1)t)(\|u_{0}\|_{H^{s}}^{2}+\|\rho_{0}\|^{2}_{H^{s-1}}+1).$$
This completes the proof of the theorem.

Given $z_{0}\in H^{s}\times H^{s-1}$ with $s\geq 2$. Theorem 2.1
ensures the existence of a maximal  $T> 0$ and a solution
                  $z=\left(\begin{array}{c}
                                      u \\
                                      \rho \\
                                    \end{array}
                                  \right)$
to (2.4) such that
$$
z=z(\cdot,z_{0})\in C([0,T); H^{s}\times H^{s-1})\cap
C^{1}([0,T);H^{s-1}\times H^{s-2}).
$$

Consider now the following initial value problem

\begin{equation}
\left\{\begin{array}{ll}q_{t}=u(t,-q)+2\gamma_{2},\ \ \ \ t\in[0,T), \\
q(0,x)=x,\ \ \ \ x\in\mathbb{R}, \end{array}\right.
\end{equation}
where $u$ denotes the first component of the solution $z$ to (2.4).
Then we have the following two useful lemmas.

Similar to the proof of Lemma 4.1 in \cite{Y1}, applying classical
results in the theory of ordinary differential equations, one can
obtain the following result on $q$ which is crucial in the proof of
blow-up scenarios.

\begin{lemma3}
 Let $u\in C([0,T); H^{s})\bigcap C^{1}([0,T); H^{s-1}), s\geq 2$. Then Eq.(3.7) has a unique solution
$q\in C^1([0,T)\times \mathbb{R};\mathbb{R})$. Moreover, the map
$q(t,\cdot)$ is an increasing diffeomorphism of $\mathbb{R}$ with
$$
q_{x}(t,x)=exp\left(-\int_{0}^{t}u_{x}(s,-q(s,x))ds\right)>0, \ \
(t,x)\in [0,T)\times \mathbb{R}.$$
\end{lemma3}

\begin{lemma3}
Let $z_{0}=\left(\begin{array}{c}
                                u_{0} \\
                                \rho_{0} \\
                              \end{array}
                            \right)
\in H^{s}\times H^{s-1}$, $s \geq 2$ and let $T>0$ be the maximal
existence time of the corresponding solution
                  $z=\left(\begin{array}{c}
                                      u \\
                                      \rho \\
                                    \end{array}
                                  \right)$
to (1.2). Then we have
\begin{equation}
\rho(t,-q(t,x))q_{x}(t,x)=\rho_{0}(-x), \ \ \ \forall \
(t,x)\in[0,T)\times \mathbb{S}.
\end{equation}
Moreover, if there exists $M> 0$ such that $u_{x}\leq M$ for all
$(t,x)\in [0,T)\times \mathbb{S}$, then
$$\|\rho(t,\cdot)\|_{L^{\infty}}\leq
e^{MT}\|\rho_{0}(\cdot)\|_{L^{\infty}}, \ \ \ \forall \ t\in[0,
T).$$
\end{lemma3}
\textbf{Proof} Differentiating the left-hand side of the equation
(3.8) with respect to $t$, and applying the relations (2.4) and
(3.7), we obtain
\begin{align*}
&\frac{d}{dt}\rho(t,-q(t,x))q_{x}(t,x)\\
=&(\rho_{t}(t,-q)-\rho_{x}(t,-q)q_{t}(t,x))q_{x}(t,x)+\rho(t,-q(t,x))q_{xt}(t,x)\\
=&(\rho_{t}-(u(t,-q)+2\gamma_{2})\rho_{x})q_{x}(t,x)-u_{x}\rho
q_{x}(t,x)\\
=&(\rho_{t}-(u+2\gamma_{2})\rho_{x}-u_{x}\rho)q_{x}(t,x)=0
\end{align*}
This proves (3.8). By Lemma 3.3, in view of (3.8) and the assumption
of the lemma, we obtain
\begin{align*}
\|\rho(t,\cdot)\|_{L^{\infty}(\mathbb{S})}&=\|\rho(t,\cdot)\|_{L^{\infty}(\mathbb{R})}\\
&=\|\rho(t,-q(t,\cdot))\|_{L^{\infty}(\mathbb{R})}\\
&=\|exp\left(\int_{0}^{t}u_{x}(s,-q(s,x))ds\right)\rho_{0}(-x)\|_{L^{\infty}(\mathbb{R})}\\
&\leq
e^{MT}\|\rho_{0}(\cdot)\|_{L^{\infty}(\mathbb{R})}=e^{MT}\|\rho_{0}(\cdot)\|_{L^{\infty}(\mathbb{S})},
\ \ \forall \ t\in[0,T).
\end{align*}

Our next result describes the precise blow-up scenario for
sufficiently regular solutions to (1.2).

\begin{theorem3}
Let $z_{0}=\left(\begin{array}{c}
                                u_{0} \\
                                \rho_{0} \\
                              \end{array}
                            \right)
\in H^{s}\times H^{s-1}$, $s>\frac{5}{2} $ be given and let T be the
maximal existence time of the corresponding solution
                  $z=\left(\begin{array}{c}
                                      u \\
                                      \rho \\
                                    \end{array}
                                  \right)
$ to (2.4) with the initial data $z_{0}$. Then the corresponding
solution blows up in finite time if and only if
$$\limsup\limits_{t\rightarrow T}\sup\limits_{x\in \mathbb{S}}\{u_{x}(t,x)\}=+\infty  \ \
\text{or} \ \ \limsup\limits_{t\rightarrow
T}\{\|\rho_{x}(t,\cdot)\|_{L^{\infty}}\}=+\infty.$$
\end{theorem3}
\textbf{Proof} By Theorem 2.1 and Sobolev's imbedding theorem it is
clear that if
$$\limsup\limits_{t\rightarrow T}\sup\limits_{x\in \mathbb{S}}\{u_{x}(t,x)\}=+\infty  \ \
\text{or} \ \ \limsup\limits_{t\rightarrow
T}\{\|\rho_{x}(t,\cdot)\|_{L^{\infty}}\}=+\infty,$$ then $T<\infty$.

Let $T<\infty$. Assume that there exists $M_{1}>0$ and $M_{2}>0$
such that
$$u_{x}(t,x)\leq M_{1}, \ \ \forall\ (t,x)\in
[0,T)\times\mathbb{S},$$ and
$$\|\rho_{x}(t,\cdot)\|_{L^{\infty}}\leq M_{2}, \ \ \ \forall \ t\in
[0,T).$$ By Lemma 3.4, we have$$\|\rho(t,\cdot)\|_{L^{\infty}}\leq
e^{M_{1}T}\|\rho_{0}\|_{L^{\infty}}, \ \ \forall\ t\in [0,T).$$ By
(2.2) and the first equation in (2.4), a direct computation implies
\begin{align}
&\frac{d}{dt}\int_{\mathbb{S}}u^{2}(t,x)dx\\
\nonumber=\
&2\int_{\mathbb{S}}u\left((u+\gamma_{1})u_{x}+\partial_{x}^{-1}(-2\mu_{0}
u-\frac{1}{2}u_{x}^{2}-\frac{1}{2}\rho^{2}+a)+h(t)\right)dx\\
\nonumber\leq \
&\int_{\mathbb{S}}u^{2}dx+\int_{\mathbb{S}}\left(\int_{0}^{x}(-2\mu_{0}
u-\frac{1}{2}u_{y}^{2}-\frac{1}{2}\rho^{2}+a)dy\right)^{2}dx+2|h(t)|\int_{\mathbb{S}}|u(t,x)|dx\\
\nonumber\leq \
&\int_{\mathbb{S}}u^{2}dx+8\mu_{0}^{2}(\int_{\mathbb{S}}|u|dx)^{2}+2\left(\int_{\mathbb{S}}
(\frac{1}{2}u_{x}^{2}+\frac{1}{2}\rho^{2}+a)dx\right)^{2}\\
&\nonumber+\max\limits_{t\in[0,T)}|h(t)|+\max\limits_{t\in[0,T)}|h(t)|\int_{\mathbb{S}}u^{2}(t,x)dx\\
\nonumber= \
&(1+8\mu_{0}^{2}+\max\limits_{t\in[0,T)}|h(t)|)\int_{\mathbb{S}}u^{2}dx
+\frac{1}{2}\left[\int_{0}^{1}(u_{0,x}^{2}+\rho_{0}^{2}+2a)dx\right]^{2}+\max\limits_{t\in[0,T)}|h(t)|
\end{align}
for $t\in(0,T)$.

Multiplying the first equation in (1.2) by $m=u_{xx}$ and
integrating by parts, we find
\begin{align}
\frac{d}{dt}\int_{\mathbb{S}}m^{2}dx= \
&-4\mu\int_{\mathbb{S}}mu_{x}dx+
4\int_{\mathbb{S}}u_{x}m^{2}dx+2\int_{\mathbb{S}}umm_{x}dx\\
&\nonumber-2\int_{\mathbb{S}}m\rho\rho_{x}dx
+2\gamma_{1}\int_{\mathbb{S}}mm_{x}dx\\
\nonumber= \
&3\int_{\mathbb{S}}u_{x}m^{2}dx-2\int_{\mathbb{S}}m\rho\rho_{x}dx\\
\nonumber\leq \
&3M_{1}\int_{\mathbb{S}}m^{2}dx+\|\rho\|_{L^{\infty}}\int_{\mathbb{S}}m^{2}+\rho_{x}^{2}dx\\
\nonumber\leq \
&(3M_{1}+\|\rho\|_{L^{\infty}})\int_{\mathbb{S}}m^{2}dx+\|\rho\|_{L^{\infty}}\int_{\mathbb{S}}\rho_{x}^{2}dx.
\end{align}

Differentiating the first equation in (1.2) with respect to $x$,
multiplying the obtained equation by $m_{x}=u_{xxx},$ integrating by
parts and using Lemma 3.4, we obtain
\begin{align}
&\frac{d}{dt}\int_{\mathbb{S}}m_{x}^{2}dx\\
\nonumber= \
&-4\mu\int_{\mathbb{S}}mm_{x}+4\int_{\mathbb{S}}m^{2}m_{x}dx+6\int_{\mathbb{S}}u_{x}m_{x}^{2}
+2\int_{\mathbb{S}}um_{xx}m_{x}\\
\nonumber&-2\int_{\mathbb{S}}\rho_{x}^{2}m_{x}-2\int_{\mathbb{S}}\rho\rho_{xx}m_{x}dx+2\gamma_{1}\int_{\mathbb{S}}m_{x}m_{xx}dx\\
\nonumber= \
&5\int_{\mathbb{S}}u_{x}m_{x}^{2}dx-2\int_{\mathbb{S}}\rho_{x}^{2}m_{x}dx-2\int_{\mathbb{S}}\rho\rho_{xx}m_{x}dx\\
\nonumber\leq \
&5M_{1}\int_{\mathbb{S}}m_{x}^{2}dx+2\|\rho_{x}\|_{L^{\infty}}^{2}\int_{\mathbb{S}}|m_{x}|dx+\|\rho\|_{L^{\infty}}\int_{\mathbb{S}}(\rho_{xx}^{2}+m_{x}^{2})dx\\
\nonumber\leq \
&5M_{1}\int_{\mathbb{S}}m_{x}^{2}dx+\|\rho\|_{L^{\infty}}\int_{\mathbb{S}}(\rho_{xx}^{2}+m_{x}^{2})dx+2\|\rho_{x}\|_{L^{\infty}}^{2}
+2\|\rho_{x}\|_{L^{\infty}}^{2}\int_{\mathbb{S}}m_{x}^{2}dx\\
\nonumber\leq \
&(5M_{1}+\|\rho\|_{L^{\infty}}+2M_{2}^{2})\int_{\mathbb{S}}m_{x}^{2}dx+
\|\rho\|_{L^{\infty}}\int_{\mathbb{S}}\rho_{xx}^{2}dx+2M_{2}^{2}.
\end{align}

Differentiating the second equation in (1.2) with respect to $x$,
multiplying the obtained equation by $\rho_{x}$ and integrating by
parts, we obtain
\begin{align}
\frac{d}{dt}\int_{\mathbb{S}}\rho_{x}^{2}dx= \
&3\int_{\mathbb{S}}u_{x}\rho_{x}^{2}dx+2\int_{\mathbb{S}}m\rho\rho_{x}dx\\
\nonumber\leq \
&3M_{1}\int_{\mathbb{S}}\rho_{x}^{2}dx+\|\rho\|_{L^{\infty}}\int_{\mathbb{S}}(m^{2}+\rho_{x}^{2})dx\\
\nonumber=\ &
(3M_{1}+\|\rho\|_{L^{\infty}})\int_{\mathbb{S}}\rho_{x}^{2}dx+\|\rho\|_{L^{\infty}}\int_{\mathbb{S}}m^{2}dx.
\end{align}

Differentiating the second equation in (1.2) with respect to $x$
twice, multiplying the obtained equation by $\rho_{xx},$ integrating
by parts and using Lemma 3.4, we obtain
\begin{align}
&\frac{d}{dt}\int_{\mathbb{S}}\rho_{xx}^{2}dx\\
\nonumber= \
&5\int_{\mathbb{S}}u_{x}\rho_{xx}^{2}dx+\int_{\mathbb{S}}u_{xxx}(2\rho\rho_{xx}-3\rho_{x}^{2})dx\\
\nonumber\leq \
&5M_{1}\int_{\mathbb{S}}\rho_{xx}^{2}dx+\int_{\mathbb{S}}m_{x}(2\rho\rho_{xx}-3\rho_{x}^{2})dx\\
\nonumber\leq \
&5M_{1}\int_{\mathbb{S}}\rho_{xx}^{2}dx+3\|\rho_{x}\|_{L^{\infty}}^{2}\int_{\mathbb{S}}|m_{x}|dx+
\|\rho\|_{L^{\infty}}\int_{\mathbb{S}}2m_{x}\rho_{xx}dx\\
\nonumber\leq \
&(5M_{1}+\|\rho\|_{L^{\infty}})\int_{\mathbb{S}}\rho_{xx}^{2}dx+(3M_{2}^{2}
+\|\rho\|_{L^{\infty}})\int_{\mathbb{S}}m_{x}^{2}dx+3M_{2}^{2}.
\end{align}

Summing (2.2) and (3.9)-(3.13), we have
\begin{align*}
&\frac{d}{dt}\int_{\mathbb{S}}(u^{2}+u_{x}^{2}+m^{2}+m_{x}^{2}+\rho^{2}+\rho_{x}^{2}+\rho_{xx}^{2})dx\\
\leq \
&K_{1}\int_{\mathbb{S}}(u^{2}+u_{x}^{2}+m^{2}+m_{x}^{2}+\rho^{2}+\rho_{x}^{2}+\rho_{xx}^{2})dx+K_{2},\\
\end{align*}
where
$$K_{1}=1+8\mu_{0}^{2}+\max\limits_{t\in[0,T)}|h(t)|+8e^{M_{1}T}\|\rho_{0}\|_{L^{\infty}}+16M_{1}+5M_{2}^{2},$$
$$K_{2}=\frac{1}{2}\left[\int_{\mathbb{S}}(u_{0,x}^{2}+\rho_{0}^{2}+2a)dx\right]^{2}+\max\limits_{t\in[0,T)}|h(t)|+5M_{2}^{2}.$$
By means of Gronwall's inequality and the above inequality, we
deduce that
\begin{align*}
&\|u(t,\cdot)\|_{H^{3}}^{2}+\|\rho(t,\cdot)\|_{H^{2}}^{2}\\
\leq \
&e^{K_{1}t}(\|u_{0}\|_{H^{3}}^{2}+\|\rho_{0}\|_{H^{2}}^{2}+\frac{K_{2}}{K_{1}}),
\ \ \ \forall \ t\in[0,T).
\end{align*}
The above inequality, Sobolev's imbedding theorem and Theorem 3.1
ensure that the solution $z$ does not blow-up in finite time. This
completes the proof of the theorem.

For initial data $z_{0}=\left(\begin{array}{c}
                                u_{0} \\
                                \rho_{0} \\
                              \end{array}
                            \right)
\in H^{2}\times H^{1}$, we have the following precise blow-up
scenario.

\begin{theorem3}
Let $z_{0}=\left(\begin{array}{c}
                                u_{0} \\
                                \rho_{0} \\
                              \end{array}
                            \right)
\in H^{2}\times H^{1}$, and let T be the maximal existence time of
the corresponding solution
                  $z=\left(\begin{array}{c}
                                      u \\
                                      \rho \\
                                    \end{array}
                                  \right)
$ to (2.4) with the initial data $z_{0}$. Then the corresponding
solution blows up in finite time if and only if
$$\limsup\limits_{t\rightarrow T}\sup\limits_{x\in \mathbb{S}}\{u_{x}(t,x)\}=+\infty.$$
\end{theorem3}
\textbf{Proof} Let $z=\left(\begin{array}{c}
                                      u \\
                                      \rho \\
                                    \end{array}
                                  \right)
$ be the solution to (2.4) with the initial data $z_{0}\in
H^{2}\times H^{1},$ and let $T$ be the maximal existence time of the
solution $z$, which is guaranteed by Theorem 2.1.

Let $T< \infty$. Assume that there exists $M_{1}>0$ such that
$$u_{x}(t,x)\leq M_{1}, \ \ \forall\ (t,x)\in
[0,T)\times\mathbb{S}.$$ By Lemma 3.4, we
have$$\|\rho(t,\cdot)\|_{L^{\infty}}\leq
e^{M_{1}T}\|\rho_{0}\|_{L^{\infty}}, \ \ \forall\ t\in [0,T).$$

Combining (2.2), (3.9)-(3.10) and (3.12), we
obtain$$\frac{d}{dt}\int_{\mathbb{S}}(u^{2}+u_{x}^{2}+m^{2}+\rho^{2}+\rho_{x}^{2})dx\leq
K_{3}\int_{\mathbb{S}}(u^{2}+u_{x}^{2}+m^{2}+\rho^{2}+\rho_{x}^{2})dx+K_{4},$$
where
$$K_{3}=1+8\mu_{0}^{2}+\max\limits_{t\in[0,T)}|h(t)|+6M_{1}+4e^{M_{1}T}\|\rho_{0}\|_{L^{\infty}},$$
$$K_{4}=\frac{1}{2}\left[\int_{0}^{1}(u_{0,x}^{2}+\rho_{0}^{2}+2a)dx\right]^{2}+\max\limits_{t\in[0,T)}|h(t)|.$$ By means of Gronwall's
inequality and the above inequality, we get
$$\|u(t,\cdot)\|_{H^{2}}^{2}+\|\rho(t,\cdot)\|_{H^{1}}^{2}\leq
e^{K_{3}t}(\|u_{0}\|_{H^{2}}^{2}+\|\rho_{0}\|_{H^{1}}^{2}+\frac{K_{4}}{K_{3}}).$$
The above inequality ensures that the solution $z$ does not blow-up
in finite time.

On the other hand, by Sobolev's imbedding theorem, we see that if
$$\limsup\limits_{t\rightarrow T}\sup\limits_{x\in \mathbb{S}}\{u_{x}(t,x)\}=+\infty,$$ then the solution will blow up in
finite time. This completes the proof of the theorem.

\begin{remark3} Note that Theorem 3.2 shows that
$$T(\|z_{0}\|_{H^{s}\times H^{s-1}})=T(\|z_{0}\|_{H^{s^{\prime}}
\times H^{s^{\prime}-1}}),\ \ \ \ \ \forall
s,s^{\prime}>\frac{5}{2},$$ while Theorem 3.3 implies
that$$T(\|z_{0}\|_{H^{s}\times H^{s-1}})\leq T(\|z_{0}\|_{H^{2}
\times H^{1}}),\ \ \ \ \ \forall s,s^{\prime}\geq 2.$$
\end{remark3}

\section{Blow-up}
\newtheorem{theorem4}{Theorem}[section]
\newtheorem{lemma4}{Lemma}[section]
\newtheorem {remark4}{Remark}[section]
\newtheorem{corollary4}{Corollary}[section]

In this section, we discuss the blow-up phenomena of the system
(1.2) and prove that there exist strong solutions to (1.2) which do
not exist globally in time.

\begin{lemma4}(\cite{constantin, fu})If $f\in H^{1}(\mathbb{S})$ is such that
$\int_{\mathbb{S}}f(x)dx=0,$ then we have
$$\max\limits_{x\in\mathbb{S}}f^{2}(x)\leq
\frac{1}{12}\int_{\mathbb{S}}f_{x}^{2}(x)dx.$$
\end{lemma4}

Note that $\int_{\mathbb{S}}(u(t,x)-\mu_{0})dx=\mu_{0}-\mu_{0}=0.$
By Lemma 4.1, we find that
$$\max\limits_{x\in\mathbb{S}}[u(t,x)-\mu_{0}]^{2}\leq
 \frac{1}{12}\int_{\mathbb{S}}u_{x}^{2}(t,x)dx\leq
 \frac{1}{12}\mu_{1}^{2}.$$ So we have
 \begin{equation}
 \|u(t,\cdot)\|_{L^{\infty}(\mathbb{S})}\leq
 |\mu_{0}|+\frac{\sqrt{3}}{6}\mu_{1}.
 \end{equation}

\begin{theorem4}Let $z_{0}=\left(
                                                     \begin{array}{c}
                                                       u_{0} \\
                                                       \rho_{0} \\
                                                     \end{array}
                                                   \right)\not\equiv
                                                   0
\in H^s\times H^{s-1}, s\geq 2,$  and T be the maximal time of the
solution $z=\left(
                                    \begin{array}{c}
                                      u \\
                                      \rho\\
                                    \end{array}
                                  \right)$
to (1.2) with the initial data $z_0$. If $\gamma_{1}=2\gamma_{2},$
$\mu_{0}=0$ and there exists a point $x_{0}\in \mathbb{S},$ such
that $\rho_{0}(-x_{0})=0$, then the corresponding solution to (1.2)
blows up in finite time.
\end{theorem4}
\textbf{Proof} Let $m(t)=u_{x}(t,-q(t,x_{0}))$, $\gamma(t)=\rho
(t,-q(t,x_{0}))$, where $q(t,x)$ is the solution of Eq.(3.7). By
Eq.(3.7) we can obtain
$$\frac{dm}{dt}=(u_{tx}-(u+\gamma_{1})u_{xx})(t,-q(t,x_{0})).$$
Evaluating the integrated representation (2.3) at $(t,-q(t,x_{0}))$
with the assumption $\mu_{0}=0$, we get
$$\frac{d}{dt}m(t)=\frac{1}{2}m(t)^{2}-\frac{1}{2}\gamma(t)^{2}+a.$$
Since $\gamma(0)=0$, we infer from Lemmas 3.3-3.4 that $\gamma(t)=0$
for all $t\in[0,T).$ Note that
$a=2\mu(u)^{2}+\frac{1}{2}\int_{\mathbb{S}}(u_{x}^{2}+\rho^{2})dx>0.$
(Indeed, if $a(t)=0,$ then $(u,\rho)=(0,0).$ This contradicts the
assumption of the theorem.) Then we have $\frac{d}{dt}m(t)\geq a >
0$. Thus, it follows that $m(t_{0})>0$ for some $t_{0}\in(0,T).$
Solving the following inequality yields
$$\frac{d}{dt}m(t)\geq \frac{1}{2}m(t)^{2}.$$ Therefore $$ 0 <
\frac{1}{m(t)}\leq \frac{1}{m(t_{0})}-\frac{1}{2}(t-t_{0}), \ \ \ \
 t\in[t_{0},T).$$ The above inequality implies that $T<
t_{0}+\frac{2}{m(t_{0})}$ and $\lim\limits_{t\rightarrow T}m(t)=
+\infty.$ In view of Theorem 3.2, this completes the proof of the
theorem.

\begin{theorem4}Let $z_{0}=\left(
                                                     \begin{array}{c}
                                                       u_{0} \\
                                                       \rho_{0} \\
                                                     \end{array}
                                                   \right)
\in H^s\times H^{s-1}, s\geq 2,$  and T be the maximal time of the
solution $z=\left(
                                    \begin{array}{c}
                                      u \\
                                      \rho\\
                                    \end{array}
                                  \right)$
to (1.2) with the initial data $z_0$. If $\gamma_{1}=2\gamma_{2},$
$\mu_{0}\neq 0,$ $|\mu_{0}|+\frac{\sqrt{3}}{6}\mu_{1}<
\frac{a}{2|\mu_{0}|}$ and there exists a point $x_{0}\in
\mathbb{S},$ such that $\rho_{0}(-x_{0})=0$, then the corresponding
solution to (1.2) blows up in finite time.
\end{theorem4}
\textbf{Proof} Let $m(t)=u_{x}(t,-q(t,x_{0}))$, $\gamma(t)=\rho
(t,-q(t,x_{0}))$, where $q(t,x)$ is the solution of Eq.(3.7). By
Eq.(3.7) we can obtain
$$\frac{dm}{dt}=(u_{tx}-(u+\gamma_{1})u_{xx})(t,-q(t,x_{0})).$$
Evaluating the integrated representation (2.3) at $(t,-q(t,x_{0}))$
we have
$$\frac{d}{dt}m(t)=\frac{1}{2}m(t)^{2}-\frac{1}{2}\gamma(t)^{2}+a-2\mu_{0}u.$$
Since $\gamma(0)=0$, we infer from Lemmas 3.3-3.4 that $\gamma(t)=0$
for all $t\in[0,T).$ In view of (4.1) and the condition
$|\mu_{0}|+\frac{\sqrt{3}}{6}\mu_{1}< \frac{a}{2|\mu_{0}|},$ we have
$a-2\mu_{0}u\geq a-2|\mu_{0}u|>0.$ Then we have
$\frac{d}{dt}m(t)\geq a-2\mu_{0}u > 0$. The left proof is the same
as that of Theorem 4.1, so we omit it here.

\section{Global Existence}
\newtheorem{theorem5}{Theorem}[section]
\newtheorem{lemma5}{Lemma}[section]
\newtheorem {remark5}{Remark}[section]
\newtheorem{corollary5}{Corollary}[section]
In this section, we will present a global existence result. Firstly,
we give two useful lemmas.

\begin{theorem5}Let $z_{0}=\left(
                                                     \begin{array}{c}
                                                       u_{0} \\
                                                       \rho_{0} \\
                                                     \end{array}
                                                   \right)
\in H^{2}\times H^{1},$ and T be the maximal time of the solution
$z=\left(
                                    \begin{array}{c}
                                      u \\
                                      \rho \\
                                    \end{array}
                                  \right)
$ to (1.2) with the initial data $z_0$. If $\gamma_{1}=2\gamma_{2},$
$\rho_{0}(x)\neq 0$ for all $x\in\mathbb{S}$, then the corresponding
solution $z$ exists globally in time.
\end{theorem5}
\textbf{Proof} By Lemma 3.3, we know that $q(t,\cdot)$ is an
increasing diffeomorphism of $\mathbb{R}$ with
$$
q_{x}(t,x)=exp\left(-\int_{0}^{t}u_{x}(s,-q(s,x))ds\right)>0, \ \
\forall \ (t,x)\in [0,T)\times \mathbb{R}.$$ Moreover,
\begin{eqnarray}
\sup_{y\in\mathbb{S}}u_{y}(t,y)=\sup_{x\in\mathbb{R}}u_{x}(t,-q(t,x)),
 \ \forall \ t\in[0,T).
\end{eqnarray}
Set $M(t,x)=u_{x}(t,-q(t,x))$ and $\alpha(t,x)=\rho(t,-q(t,x))$ for
$t\in [0,T)$ and $x\in\mathbb{R}$. By $\gamma_{1}=2\gamma_{2},$
(1.2) and Eq.(3.7), we have
\begin{eqnarray}
\frac{\partial M}{\partial
t}=(u_{tx}-(u+\gamma_{1})u_{xx})(t,-q(t,x))\ \ \text{and} \ \
\frac{\partial \alpha}{\partial t}=\alpha M.
\end{eqnarray}
Evaluating (2.3) at $(t,-q(t,x))$ we get
$$
\partial_{t}M(t,x)=\frac{1}{2}M(t,x)^{2}-\frac{1}{2}\alpha(t,x)^{2}+a-2\mu_{0} u(t,-q(t,x)).
$$
Write $f(t,x)=a-2\mu_{0} u(t,-q(t,x)).$ By (4.1) we have
\begin{align*}
|f(t,x)|\leq a+2|\mu_{0}|\|u\|_{L^{\infty}}&\leq
a+2|\mu_{0}|(|\mu_{0}|+\frac{\sqrt{3}}{6}\mu_{1})\\
&=4\mu_{0}^{2}+\frac{1}{2}\mu_{1}^{2}+\frac{\sqrt{3}}{3}
|\mu_{0}|\mu_{1}
\end{align*} and
\begin{equation}
\partial_{t}M(t,x)=\frac{1}{2}M(t,x)^{2}-\frac{1}{2}\alpha(t,x)^{2}+f(t,x).
\end{equation}

 By Lemmas 3.3-3.4, we know that $\alpha(t,x)$ has the same sign
with $\alpha(0,x)=\rho_{0}(-x)$ for every $x\in \mathbb{R}$.
Moreover, there is a constant $\beta>0$ such that
$\inf\limits_{x\in\mathbb{R}}|\alpha(0,x)|=\inf\limits_{x\in\mathbb{S}}|\rho_{0}(-x)|\geq\beta>0$
since $\rho_{0}(x)\neq 0$ for all $x\in\mathbb{S}$ and $\mathbb{S}$
is a compact set. Thus,
$$\alpha(t,x)\alpha(0,x)>0, \ \ \ \forall x\in \mathbb{R}.$$Next, we
consider the following Lyapunov function first introduced in
\cite{C-I}.
\begin{equation}
w(t,x)=\alpha(t,x)\alpha(0,x)+\frac{\alpha(0,x)}{\alpha(t,x)}(1+M^{2}),
\ \ \ (t,x)\in[0,T)\times\mathbb{R}.
\end{equation}
By Sobolev's imbedding theorem, we have
\begin{align}
0&<w(0,x)=\alpha(0,x)^{2}+1+M(0,x)^{2}\\
&\nonumber=\rho_{0}(x)^{2}+1+u_{0,x}(x)^{2}\\
&\nonumber\leq1+\max\limits_{x\in\mathbb{S}}(\rho_{0}(x)^{2}+u_{0,x}(x)^{2}):=C_{1}.
\end{align}
 Differentiating
(5.4) with respect to $t$ and using (5.2)-(5.3), we obtain

\begin{align}
\nonumber \frac{\partial w}{\partial
t}(t,x)&=\frac{\alpha(0,x)}{\alpha(t,x)}M(t,x)(2f-1)\\
&\nonumber
\leq|f-\frac{1}{2}|\frac{\alpha(0,x)}{\alpha(t,x)}(1+M^{2})\\
&\nonumber\leq(4\mu_{0}^{2}+\frac{1}{2}\mu_{1}^{2}+\frac{\sqrt{3}}{3}
|\mu_{0}|\mu_{1}+\frac{1}{2})w(t,x).
\end{align}
By Gronwall's inequality, the above inequality and (5.5), we have
$$
w(t,x)\leq
w(0,x)e^{(4\mu_{0}^{2}+\frac{1}{2}\mu_{1}^{2}+\frac{\sqrt{3}}{3}
|\mu_{0}|\mu_{1}+\frac{1}{2})t}\leq
C_{1}e^{(4\mu_{0}^{2}+\frac{1}{2}\mu_{1}^{2}+\frac{\sqrt{3}}{3}
|\mu_{0}|\mu_{1}+\frac{1}{2})t}
$$
for all $(t,x)\in [0,T)\times\mathbb{R}.$ On the other
hand,$$w(t,x)\geq2\sqrt{\alpha^{2}(0,x)(1+M^{2})}\geq2\beta|M(t,x)|,
\ \ \ \forall \ \ (t,x)\in [0,T)\times\mathbb{R}.$$ Thus,
$$
|M(t,x)|\leq \frac{1}{2\beta}w(t,x)\leq
\frac{1}{2\beta}C_{1}e^{(4\mu_{0}^{2}+\frac{1}{2}\mu_{1}^{2}+\frac{\sqrt{3}}{3}
|\mu_{0}|\mu_{1}+\frac{1}{2})t}
$$
for all $(t,x)\in [0,T)\times\mathbb{R}$. Then by (5.1) and the
above inequality, we have
$$\limsup\limits_{t\rightarrow
T}\sup\limits_{y\in\mathbb{S}}u_{y}(t,y)=\limsup\limits_{t\rightarrow
T}\sup\limits_{x\in\mathbb{R}}u_{x}(t,-q(t,x))
\leq\frac{1}{2\beta}C_{1}e^{(4\mu_{0}^{2}+\frac{1}{2}\mu_{1}^{2}+\frac{\sqrt{3}}{3}
|\mu_{0}|\mu_{1}+\frac{1}{2})t}.
$$
This completes the proof by using Theorem 3.3.

\bigskip

\noindent\textbf{Acknowledgments} This work was partially supported
by NNSFC (No. 10971235), RFDP (No. 200805580014), NCET-08-0579 and
the key project of Sun Yat-sen University.


\begin{thebibliography}{99}
\small



\bibitem{B1}R. Beals, D. Sattinger and J. Szmigielski, Inverse scattering solutions of the Hunter-Saxton
equations, {\it Appl. Anal.}, {\bf78} (2001), 255-269.

\bibitem{R-D}R. Camassa and D. Holm, An integrable shallow water equation with peaked solitons, {\it Phys.
Rev. Lett.}, {\bf 71} (1993), 1661-1664.

\bibitem{constantin}
A. Constantin, On the Blow-up of solutions of a periodic shallow
water equation, {\it J. Nonlinear Sci.}, {\bf 10} (2000), 391-399.

\bibitem{C-I}
A. Constantin and R. I. Ivanov, On an integrable two-component
Camass-Holm shallow water system, {\it Phys. Lett. A}, {\bf 372}
(2008), 7129--7132.

\bibitem{C-K}A. Constantin and B. Kolev, On the geometric approach to the motion of inertial mechanical
systems, {\it J. Phys. A}, {\bf 35} (2002), R51-R79.

\bibitem{C-M}A. Constantin and H. P. McKean, A shallow water equation on the circle, {\it Comm. Pure
Appl. Math.}, {\bf 52} (1999), 949-982.

\bibitem{HHD}H. H. Dai and M. Pavlov, Transformations for the Camassa-Holm equation, its high-frequency
limit and the Sinh-Gordon equation, {\it J. P. Soc. Japan}, {\bf 67}
(1998), 3655-3657.


\bibitem{E-L-Y}
J. Escher, O. Lechtenfeld and Z. Yin, Well-posedness and blow-up
phenomena for the 2-component Camassa-Holm equation, {\it Discrete
Contin. Dyn. Syst.}, {\bf 19} (2007), 493-513.

\bibitem{F-F}A. Fokas and B. Fuchssteiner, Symplectic structures, their B¡§acklund transformations and
hereditary symmetries, {\it Phys. D}, {\bf 4} (1981/1982), 47-66.

\bibitem{fu} Y. Fu, Y. Liu and C. Qu, On the blow up structure for the generalized periodic Camassa-Holm and
Degasperis-Procesi equations, arXiv:1009.2466.

\bibitem{J-R} J. K. Hunter and R. Saxton, Dynamics of director fields, {\it SIAM J. Appl. Math.}, {\bf51} (1991),
1498-1521.

\bibitem{J-Y} J. K. Hunter and Y. Zheng, On a completely integrable nonlinear
hyperbolic variational equation, {\it Phys. D}, {\bf 79} (1994),
361-386.

\bibitem{J-H}R. S. Johnson, Camassa-Holm, Korteweg-de Vries and related models for water waves, {\it J.
Fluid. Mech.}, {\bf 455} (2002), 63-82.

\bibitem{Kato 1}
T. Kato, Quasi-linear equations of evolution, with applications to
partial differential equations, in "Spectral Theory and Differential
Equations", Lecture Notes in Math., Vol. 448, Springer Verlag,
Berlin, (1975), 25--70.

\bibitem{Kato 3}
T. Kato and G. Ponce, Commutator estimates and Navier-Stokes
equations, {\it Comm. Pure Appl. Math.}, {\bf 41} (1988), 203--208.

\bibitem{k-l-m}
B. Khesin, J. Lenells and G. Misiolek, Generalized Hunter-Saxton
equation and the geometry of the group of circle diffeomorphisms,
{\it Math. Ann.}, {\bf 342} (2008), 617-656.

\bibitem{l-m-t}
J. Lenells, G. Misiolek and F. Ti$\breve{g}$lay, Integrable
evolution equations on spaces of tensor densities and their peakon
solutions, {\it Commun. Math. Phys.}, {\bf 299} (2010), 129-161.

\bibitem{l-y}
J. Liu and Z. Yin, Blow-up phenomena and global existence for a
periodic two-component Hunter-Saxton system, arXiv:1012.5448.

\bibitem{P-P} P. Olver and P. Rosenau, Tri-Hamiltonian duality between solitons and solitary wave solutions
having compact support, {\it Phys. Rev. E(3)}, {\bf53} (1996),
1900-1906.

\bibitem{MW} M. Wunsch, On the Hunter-Saxton system, {\it Discrete
Contin. Dyn. Syst. B}, {\bf 12} (2009), 647-656.

\bibitem{Y1} Z. Yin, Global existence for a new periodic integrable
equation, {\it J. Math. Anal. Appl.}, {\bf 283}, (2003), 129-139.

\bibitem{Y} Z. Yin, On the structure of solutions to the periodic Hunter-Saxton equation, {\it SIAM J. Math.
Anal.}, {\bf 36} (2004), 272-283.

\bibitem{d-z} D. Zuo, A 2-component $\mu$-Hunter-Saxton Equation,
{\it Inverse Problems}, {\bf 26} (2010), 085003 (9pp).

\end{thebibliography}
\end{document}